\documentclass[preprint]{elsarticle}
\usepackage{amsfonts}
\usepackage{amsmath}
\usepackage[all]{xy}
\usepackage{tikz}
\usetikzlibrary{trees}
\usetikzlibrary{arrows,shapes,snakes,automata,backgrounds,petri}

\newtheorem{thm}{Theorem}
\newtheorem{lem}[thm]{Lemma}
\newtheorem{cor}[thm]{Corollary}
\newproof{pf}{Proof}
\newproof{pot1}{Proof of Theorem \ref{mainimp}}
\newproof{pot2}{Proof of Theorem \ref{mainIgraph}}
\newproof{pot3}{Proof of Theorem \ref{mainIgrapha}}

\begin{document}
\begin{frontmatter}

\title{Graph-theoretic conditions for injectivity of functions on rectangular domains}

\author[ref1]{Murad Banaji\corref{cor1}}

\address[ref1]{Department of Mathematics, University College London, Gower Street, London~WC1E~6BT, UK.}
\cortext[cor1]{m.banaji@ucl.ac.uk. Research funded by EPSRC grant EP/C500210/1}

\begin{abstract}
This paper presents sufficient graph-theoretic conditions for injectivity of collections of differentiable functions on rectangular subsets of $\mathbb{R}^n$. The results have implications for the possibility of multiple fixed points of maps and flows. Well-known results on systems with signed Jacobians are shown to be easy corollaries of more general results presented here. 
\end{abstract}


\end{frontmatter}

\section{Introduction}

Finding conditions for global injectivity of functions satisfying only local, structural, conditions is of both theoretical and practical importance. A variety of such conditions, spectral and otherwise, have been found (\cite{fernandes,smyth} for example). The explicit aim of such work is often to find restrictions on the Jacobian which guarantee injectivity of certain classes of functions (e.g. maps $F:\mathbb{R}^n \to \mathbb{R}^n$ which are polynomial, analytic, or $C^1$). However there are also close links between some of this work and questions of global stability in dynamical systems. Often the questions posed are highly nontrivial, and a number of open conjectures exist in this area, for example, the ``Chamberland conjecture'', open at the time of writing, that a $C^1$ map $F:\mathbb{R}^n \to \mathbb{R}^n$ with Jacobian $DF$ whose spectrum is bounded away from the origin is injective \cite{chamberland}. This conjecture, if proved, would imply a number of other injectivity results.

Apart from the theoretical interest, injectivity of functions is important in a variety of applications: in particular for exploring the possibility of multiple fixed points of maps or flows. Although spectral conditions seem to have been of most theoretical interest, one particular class of injective functions -- differentiable functions on a rectangular domain with $P$ matrix Jacobians \cite{gale} (notions to be defined below) -- has proved of relevance in several practical contexts \cite{soule,banajiSIAM,banajicraciun2}. While $P$ matrices have spectra disjoint from a region of the complex plane \cite{kellogg}, they are not defined by their spectra, although some results in \cite{smyth} can be interpreted as generalisations of certain $P$ matrix results which impose no global coordinates. \cite{parthasarathy} provides a useful summary of results in this area, along with generalisations and applications. 

The primary aim of this paper is to present graph-theoretic corollaries of the $P$ matrix results, and construct links between two strands of theory: work on injectivity of functions with signed Jacobian (e.g. \cite{soule,gouze98,kaufman}) and work on injectivity of more general functions in \cite{banajicraciun2} which extended earlier work in \cite{banajiSIAM,craciun,craciun1,banajicraciun}. These approaches use different generalisations of graphs -- the former use ``interaction graphs'' (here abbreviated to ``I-graphs''), while the latter use variants of the so-called ``SR graph'', originally defined for dynamical systems arising from systems of chemical reactions \cite{craciun1}. Both I-graphs and a directed variant of SR graphs, termed DSR graphs, will be defined below. The relevant previous results on I-graphs and SR graphs can be summarised as: 
\begin{itemize}
\item[A.] Injectivity of certain functions can be deduced by constructing the I-graphs associated with these functions, and confirming that these contain no positive cycles.
\item[B.] Injectivity of certain functions can be deduced by constructing DSR graphs associated with these functions, and checking conditions on cycles in these graphs.
\end{itemize}

The I-graph and DSR graph are both derived from Jacobians of the functions, and the results are closely related to the question of when these Jacobians are $P$ matrices or in some closely related class. In \cite{kaufman}, Kaufmann {\em et al} commented that the approaches are unrelated. Here it will be shown that, on the contrary, a number of I-graph results are corollaries of DSR graph results. The main results of this paper are:
\begin{enumerate}
\item A strengthening of DSR graph results on injectivity in \cite{banajicraciun2}: the key idea is to include certain ``nondegeneracy'' conditions on DSR graphs allowing one to enlarge the set of functions to which statement B applies. This enlargement is carried out in Theorem~\ref{nond}.
\item Theorem~\ref{mainimp}, which states that any conclusions about injectivity that can be drawn from the absence or presence of positive (resp. negative) cycles in I-graphs, are a subset of results which can be derived from DSR graphs. In other words, the functions to which statement A applies are a proper subset of those to which statement B applies.
\end{enumerate}

\section{Basic notions}

A function $f:X \to Y$ is {\bf injective} on $X$ if $f(x_1) = f(x_2)$ implies $x_1 = x_2$ for $x_1, x_2 \in X$. If $X \subset \mathbb{R}^n$ and $f:X \to \mathbb{R}^n$ is injective, then the differential equation $\dot x = f(x)$ can have no more than one equilibrium in $X$, and similarly, the map $g(x) = x + f(x)$ can have no more than one fixed point in $X$. 

A {\bf rectangular subset} of $\mathbb{R}^n$ is the product of $n$ intervals. These intervals may be closed or nonclosed, bounded or unbounded. A {\bf generalised graph} will be used to refer to a graph or multigraph, possibly directed, and possibly with additional structures including signs and labels on its edges.

{\bf Notation.} From here on, the following notation will be used:
\begin{itemize}
\item $X$ is an arbitrary rectangular subset of $\mathbb{R}^n$.
\item $\mathcal{D}_d^{+}(X)$ is the set of all differentiable\footnote{Some previous work \cite{banajicraciun2} assumed, for convenience, that all functions in question were $C^1$. It should be noted that all the results used or presented here require only differentiability and not continuous differentiability.}, diagonal functions on $X$ with range $\mathbb{R}^n$ and having positive slope, that is all $q \in \mathcal{D}_d^{+}(X)$ are of the form $q = [q_1(x_1), \ldots, q_n(x_n)]^T$ with $\frac{\partial q_i}{\partial x_i} > 0$ everywhere on $X$. 
\item $f:X \to \mathbb{R}^n$ is a differentiable function. $\mathcal{D}(X)$ is the set of all such functions. $\mathcal{F} \subset \mathcal{D}(X)$ is some collection of such functions. Define $f^{-} = \{f - q\,|\, q \in \mathcal{D}_d^{+}(X)\}$, $f^{+} = \{f + q\,|\, q \in \mathcal{D}_d^{+}(X)\}$, $\mathcal{F}^{-} = \{f - q\,|\, f \in \mathcal{F}, q \in \mathcal{D}_d^{+}(X)\}$, and $\mathcal{F}^{+} = \{f + q\,|\, f \in \mathcal{F}, q \in \mathcal{D}_d^{+}(X)\}$. 
\end{itemize}

{\bf Key goals.} A collection $\mathcal{F} \subset \mathcal{D}(X)$ will be termed injective on $X$ if $f$ is injective on $X$ for each $f \in \mathcal{F}$. Given $\mathcal{F} \subset \mathcal{D}(X)$, generalised graphs associated with $\mathcal{F}$ will be examined to make claims about injectivity of $\mathcal{F}^{-}$, which can in some cases be extended to claims about injectivity of $\mathcal{F}$. All the results have dual versions: for each claim about $\mathcal{F}^{-}$, there is a corresponding claim about $\mathcal{F}^{+}$. These dual results are collected in~\ref{dualapp}. They follow naturally from the main results and proofs are omitted.

{\bf Matrices: notation and definitions.} Let $M$ be an $n \times m$ matrix, and $\gamma \subset \{1, \ldots, n\}$, $\delta \subset \{1, \ldots, m\}$ nonempty sets. $M(\gamma|\delta)$ is the submatrix of $M$ with rows indexed by $\gamma$ and columns indexed by $\delta$. A principal submatrix of $M$ is a submatrix of the form $M(\gamma|\gamma)$. If $|\gamma| = |\delta|$, then $M[\gamma|\delta]$ means $\mathrm{det}(M(\gamma|\delta))$. {\bf Principal minors} are determinants of principal submatrices: $M[\gamma]$ is shorthand for $M[\gamma|\gamma]$. {\bf $P$ matrices} are square matrices all of whose principal minors are positive. They are by definition nonsingular. $M$ determines the {\bf qualitative class} $\mathcal{Q}(M)$ \cite{brualdi} of all matrices with the same sign pattern as $M$. Explicitly, $\mathcal{Q}(M)$ consists of all $n \times m$ matrices $X$ satisfying $M_{ij}X_{ij} > 0$ when $M_{ij} \not=0$, and $X_{ij} = 0$ when $M_{ij} = 0$. A square matrix $M$ is {\bf sign nonsingular} if all matrices in $\mathcal{Q}(M)$ are nonsingular. 



\section{I-graphs: construction and results}

For maximum generality, an ``I-graph'' is defined to be a directed multigraph on $n$ vertices where each edge has a sign ($+1$ or $-1$). Any $n \times n$ matrix $J$, is associated with an I-graph on $n$ vertices, $H_J$, in a way which is well known: if $J_{ij} > 0$, then there is a positive directed edge in $H_J$ from vertex $j$ to vertex $i$; if  $J_{ij} < 0$, then there is a negative directed edge from vertex $j$ to vertex $i$, and if $J_{ij} = 0$, then there is no directed edge in from vertex $j$ to vertex $i$. Note that diagonal entries in $J$ correspond to self-edges in $H_J$.

Any set of $n \times n$ matrices $\mathcal{J}$ is also associated with an I-graph on $n$ vertices, $H_\mathcal{J}$, constructed by a ``superposition'' of $H_J$, for $J \in \mathcal{J}$. More precisely, $H_{\mathcal{J}}$ has a positive (resp. negative) directed edge from vertex $j$ to vertex $i$ if and only if there exists $J \in \mathcal{J}$ such that $H_J$ has a positive (resp. negative) directed edge from vertex $j$ to vertex $i$. $H_{\mathcal{J}}$ can have up to two directed edges from vertex $j$ to vertex $i$, one positive and one negative. 

Directed paths and directed cycles in I-graphs are defined in the natural way. The {\it sign} of a cycle is the product of signs of edges in the cycle. Thus a cycle is positive if it contains an even number of negative edges.

Consider a function $f \in \mathcal{D}(X)$ with Jacobian $Df(x)$. Let $\mathcal{J}_f = \{Df(x)\,|\, x \in X\}$, and define $H_f \equiv H_{\mathcal{J}_f}$. Given any $\mathcal{F} \subset \mathcal{D}(X)$, define $\mathcal{J}_{\mathcal{F}} = \{Df(x)\,|\, f \in \mathcal{F}, x \in X\}$, and $H_{\mathcal{F}} \equiv H_{\mathcal{J}_{\mathcal{F}}}$. The following results hold.

\begin{thm}
\label{mainIgraph}
Given $f \in \mathcal{D}(X)$, suppose there exists some $q \in \mathcal{D}_d^{+}(X)$, and $a, b \in X$ ($a \not = b$) such that $f(a)-q(a) = f(b) - q(b)$. Then there exists $c \in X$ such that $H_{Df(c)}$ contains a positive cycle, and thus $H_{f}$ contains a positive cycle.
\end{thm}
The following theorem is one example of how, with additional assumptions, injectivity can be extended from $f^{-}$ to $f$. 
\begin{thm}
\label{mainIgrapha}
Given $f \in \mathcal{D}(X)$ such that $Df$ has negative diagonal elements (i.e. $\frac{\partial f_i}{\partial x_i} < 0$ at each point in $X$), suppose there exist $a, b \in X$ ($a \not = b$) such that $f(a)= f(b)$. Then there exists $c \in X$ such that $H_{Df(c)}$ contains a positive cycle, and thus $H_{f}$ contains a positive cycle.
\end{thm}
The following corollary follows immediately from the previous theorems. 
\begin{cor}
\label{corIgraph}
For some $\mathcal{F} \subset \mathcal{D}(X)$, assume that $H_\mathcal{F}$ contains no positive cycles. \\
{\bf 1.} Then $\mathcal{F}^{-}$ is injective. \\
{\bf 2.} Assume in addition that $\frac{\partial f_i}{\partial x_i} < 0$ for each $f \in \mathcal{F}$ and each point in $X$. Then $\mathcal{F}^{-} \cup \mathcal{F}$ is injective.
\end{cor}
\begin{pf}
The claims follow from Theorems~\ref{mainIgraph}~and~\ref{mainIgrapha} by noting that if $H_\mathcal{F}$ contains no positive cycles, then, for each $f \in \mathcal{F}$, $H_f$ contains no positive cycles. \qquad \qed
\end{pf}

All of these results are well known and stated in a variety of slightly different forms in the literature \cite{gouze98,soule}. Theorems~\ref{mainIgraph}~and~\ref{mainIgrapha} will be proved later as corollaries of stronger results on another generalised graph, termed a DSR graph.


\section{DSR graphs: construction}

Let $n, m \in \mathbb{N}$. Let $(A, B)$ be an ordered pair of real $n \times m$ matrices. Associated with $(A, B)$ is a generalised graph, $G_{A, B}$, termed a DSR graph. Before defining $G_{A, B}$ we note its properties:
\begin{enumerate}
\item $G_{A, B}$ is bipartite with two vertex-sets: a set of $n$ vertices termed ``S-vertices''; and a set of $m$ vertices termed ``R-vertices'' (such a graph will be referred to as an $n \times m$ DSR graph). No edges can exist between two S-vertices, or between two R-vertices.
\item $G_{A, B}$ is a multigraph with up to two edges between a pair of vertices. 
\item Each edge has up to two ``directions'': S-to-R direction, R-to-S direction or both, in which case we term it an undirected edge.
\item Each edge has a sign. If two edges exist between a pair of vertices, then one is positive and one is negative.
\item Each edge has an edge-label $l$ satisfying $0< l \leq \infty$. (The label $\infty$ is used only to indicate the lack of a label $0< l < \infty$.) $\mathrm{val}(e)$ will refer to the edge-label of edge $e$. 
\end{enumerate}

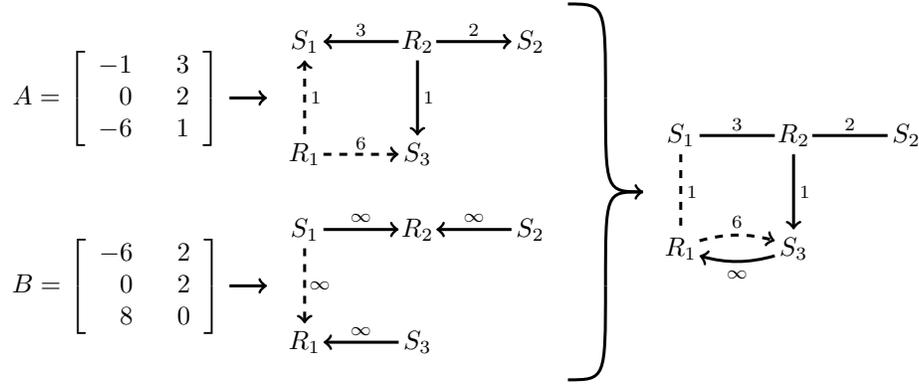
\begin{figure}[h]
\begin{center}
\begin{tikzpicture}[domain=0:4,scale=0.5]

\draw[very thin,color=black!0] (-7, 0) grid (16, 10);

\node at (-4,7.5) {$A = \left[\begin{array}{rr}-1\,\, & \,\,3\\0\,\, & \,\,2\\-6\,\, & \,\,1 \end{array}\right]$};

\draw[->, line width=0.04cm] (-1, 7.5) -- (0,7.5);

\node at (1,9) {$S_1$};
\node at (4,9) {$R_2$};
\node at (7,9) {$S_2$};
\node at (1,6) {$R_1$};
\node at (4,6) {$S_3$};

\draw[<-, line width=0.04cm] (1.5, 9) -- (3.5,9);

\draw[->, dashed, line width=0.04cm] (1.5, 6) -- (3.5,6);
\draw[->, line width=0.04cm] (4.5, 9) -- (6.5,9);

\draw[<-, dashed, line width=0.04cm] (1, 8.5) -- (1,6.5);
\draw[->, line width=0.04cm] (4, 8.5) -- (4,6.5);

\node at (2.5, 9.3) {$\scriptstyle{3}$};
\node at (2.5, 6.3) {$\scriptstyle{6}$};

\node at (1.3, 7.5) {$\scriptstyle{1}$};
\node at (4.3, 7.5) {$\scriptstyle{1}$};

\node at (5.5, 9.3) {$\scriptstyle{2}$};

\node at (-4,2.5) {$B = \left[\begin{array}{rr}-6\,\, & \,\,2\\0\,\, & \,\,2\\8\,\, & \,\,0 \end{array}\right]$};

\draw[->, line width=0.04cm] (-1, 2.5) -- (0,2.5);

\node at (1,4) {$S_1$};
\node at (4,4) {$R_2$};
\node at (7,4) {$S_2$};
\node at (1,1) {$R_1$};
\node at (4,1) {$S_3$};

\draw[->, line width=0.04cm] (1.5, 4) -- (3.5,4);

\draw[<-, line width=0.04cm] (1.5, 1) -- (3.5,1);
\draw[<-, line width=0.04cm] (4.5, 4) -- (6.5,4);

\draw[->, dashed, line width=0.04cm] (1, 3.5) -- (1,1.5);

\node at (2.5, 4.3) {$\scriptstyle{\infty}$};
\node at (2.5, 1.3) {$\scriptstyle{\infty}$};

\node at (1.4, 2.5) {$\scriptstyle{\infty}$};

\node at (5.5, 4.3) {$\scriptstyle{\infty}$};

\node at (11,6.5) {$S_1$};
\node at (14,6.5) {$R_2$};
\node at (17,6.5) {$S_2$};
\node at (11,3.5) {$R_1$};
\node at (14,3.5) {$S_3$};

\draw[-, line width=0.04cm] (11.5, 6.5) -- (13.5,6.5);

\draw[->, dashed, line width=0.04cm] (11.5, 3.7) .. controls (12.2,3.9) and (12.8, 3.9) .. (13.5,3.7);
\draw[<-, line width=0.04cm] (11.5, 3.3) .. controls (12.2,3.1) and (12.8, 3.1) .. (13.5,3.3);

\draw[-, line width=0.04cm] (14.5, 6.5) -- (16.5,6.5);

\draw[-, dashed, line width=0.04cm] (11, 6) -- (11,4);
\draw[->, line width=0.04cm] (14, 6) -- (14,4);

\node at (12.5, 6.8) {$\scriptstyle{3}$};
\node at (12.5, 4.2) {$\scriptstyle{6}$};
\node at (12.5, 2.8) {$\scriptstyle{\infty}$};

\node at (11.3, 5) {$\scriptstyle{1}$};
\node at (14.3, 5) {$\scriptstyle{1}$};

\node at (15.5, 6.8) {$\scriptstyle{2}$};

\draw[-, line width=0.04cm] (8, 0) .. controls (9,0) and (9, 0) .. (9,2.5);
\draw[->, line width=0.04cm] (9, 2.5) .. controls (9,5) and (9, 5) .. (10,5);
\draw[-, line width=0.04cm] (8, 10) .. controls (9,10) and (9, 10) .. (9,7.5);
\draw[-, line width=0.04cm] (9, 7.5) .. controls (9,5) and (9, 5) .. (10,5);

\end{tikzpicture}
\end{center}
\caption{\label{DSRbasic}Construction of a DSR graph $G_{A, B}$ from a pair of matrices $(A, B)$. Negative edges are represented as dashed lines while positive edges are bold lines, a convention which will be followed throughout. $A$ gives rise to a subgraph in which all edges have R-to-S direction, while $B$ gives a subgraph in which all edges have S-to-R direction, and edge-labels are $\infty$. A superposition of these two objects gives the DSR graph $G_{A, B}$ to the right. Note that two oppositely directed edges of the same sign in the subgraphs (e.g. between $S_1$ and $R_2$) combine to give a single undirected edge, while two oppositely directed edges with different signs (e.g. between $R_1$ and $S_3$) combine to give a pair of edges.}
\end{figure}

Since an $n \times m$ DSR graph is associated with $n \times m$ matrices, it makes sense to refer to ``S-vertex $i$'' as the S-vertex corresponding to row $i$, and ``R-vertex $j$'' as the R-vertex corresponding to column $j$. If $A_{ij} \not = 0$ and $B_{ij} = 0$, there is a single edge between R-vertex $j$ and S-vertex $i$, with R-to-S direction, the sign of $A_{ij}$, and label $|A_{ij}|$. If $B_{ij} \not = 0$ and $A_{ij} = 0$, there is a single edge between S-vertex $i$ and R-vertex $j$ with S-to-R direction with the sign of $B_{ij}$ and edge-label $\infty$. If $A_{ij}B_{ij} > 0$, then there is a single undirected edge between S-vertex $i$ and R-vertex $j$ with the sign of $A_{ij}$ and label $|A_{ij}|$. If $A_{ij}B_{ij} < 0$, then there are two edges between S-vertex $i$ and R-vertex $j$, one with R-to-S direction, the sign of $A_{ij}$, and label $|A_{ij}|$, and one with S-to-R direction, the sign of $B_{ij}$ and edge-label $\infty$. More intuition and detail are presented in \cite{banajicraciun2}. Figure~\ref{DSRbasic} provides an example of the construction.

{\bf Notation.} An edge in a DSR graph $G$ between S-vertex $i$ and R-vertex $j$ will be termed $g_{ij}$. If $g_{ij}$ has S-to-R direction it can be represented as $\overrightarrow{g}_{ij}$. Similarly if $g_{ij}$ has R-to-S direction it can be represented as $\overleftarrow{g}_{ij}$. If it is known to have both directions it can be written $\overline{g}_{ij}$. Note that referring to an edge as $\overrightarrow{g}_{ij}$ tells us that $g_{ij}$ has S-to-R direction, but does not rule out that it may also have R-to-S direction.

{\bf DSR graphs for matrix-sets.} As with I-graphs, a DSR graph can be associated with a set of matrix-pairs by taking the superposition of the DSR graphs associated with each pair. Given two sets of $n \times m$ matrices, $\mathcal{A}$ and $\mathcal{B}$, the $n \times m$ DSR graph $G_{\mathcal{A}, \mathcal{B}}$ is defined by the following requirements:
\begin{itemize}
\item If for some $A \in\mathcal{A}$, $A_{ij}\not = 0$, then $G_{\mathcal{A}, \mathcal{B}}$ contains an edge $g_{ij}$ with R-to-S direction and the sign of $A_{ij}$. Similarly if for some $B \in \mathcal{B}$, $B_{ij}\not = 0$, then $G_{\mathcal{A}, \mathcal{B}}$ contains an edge $g_{ij}$ with S-to-R direction and the sign of $B_{ij}$.
\item $G_{\mathcal{A}, \mathcal{B}}$ contains a positive (resp. negative) edge $\overleftarrow{g}_{ij}$ or $\overline{g}_{ij}$ with edge-label $0 < l <  \infty$ if and only if $A_{ij} = l$ (resp. $A_{ij} = -l$) for each $A \in \mathcal{A}$. Otherwise the edge $\overleftarrow{g}_{ij}$ or $\overline{g}_{ij}$ (if it exists) has edge-label $\infty$. An edge $\overrightarrow{g}_{ij}$ with {\em only} S-to-R direction must have edge-label $\infty$.
\end{itemize}

{\bf Properties of cycles.} Since all edges in a DSR graph are signed, all paths, and hence all cycles, have a sign defined as the product of signs of edges in the path. Define the {\bf parity} of any path $E$ of even length to be
\[
P(E) := (-1)^{|E|/2}\mathrm{sign}(E).
\]
$E$ is even if $P(E) = 1$, and odd otherwise. All cycles are paths of
even length and hence either even or odd. Even cycles are termed {\bf
  e-cycles}, while odd cycles are termed {\bf o-cycles}. A cycle $C = [e_1, e_2, \ldots, e_{2r}]$ (i.e. such that edges $e_i$ and $e_{(i \mod 2r) + 1}$ are adjacent for each $i =1,\ldots,2r$) is an {\bf s-cycle} if each edge in $C$ has a finite edge-label,  and moreover
\[
\prod_{i = 1}^{r}\mathrm{val}(e_{2i-1}) = \prod_{i = 1}^{r}\mathrm{val}(e_{2i}).
\]

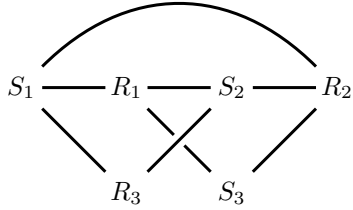
\begin{figure}[h]
\begin{center}
\begin{tikzpicture}[domain=0:4,scale=1.4]


\node at (0, 1) {$S_1$};
\node at (1, 1) {$R_1$};
\node at (2, 1) {$S_2$};
\node at (3, 1) {$R_2$};

\node at (1, 0) {$R_3$};
\node at (2, 0) {$S_3$};

\draw[-, line width=0.04cm] (0.2, 1) -- (0.8,1);
\draw[-, line width=0.04cm] (1.2, 1) -- (1.8,1);
\draw[-, line width=0.04cm] (2.2, 1) -- (2.8,1);

\draw[-, line width=0.04cm] (0.2, 0.8) -- (0.8,0.2);
\draw[-, line width=0.04cm] (2.2, 0.2) -- (2.8,0.8);
\draw[-, line width=0.04cm] (1.2, 0.2) -- (1.8,0.8);
\draw[-, line width=0.04cm] (1.2, 0.8) -- (1.45,0.55);
\draw[-, line width=0.04cm] (1.55, 0.45) -- (1.8, 0.2);

\draw[-, line width=0.04cm] (0.2, 1.2) .. controls (1,2) and (2, 2) .. (2.8, 1.2);

\end{tikzpicture}
\end{center}
\caption{\label{orient} A $3 \times 3$ DSR graph with edge-labels omitted. The cycles $C = S_1\!-\!R_1\!-\!S_3\!-\!R_2\!-\!S_2\!-\!R_3$ and $D = S_1\!-\!R_1\!-\!S_2\!-\!R_2$ have no compatible orientation, i.e. there is no choice of orientation for $C$ and $D$ such that both $S_1\!-\!R_1$ and $S_2\!-\!R_2$ have the same $C$-orientation and $D$-orientation. Thus $C$ and $D$ do not have S-to-R intersection.}
\end{figure}

{\bf Orientation of cycles.} If a cycle $C$ in a DSR graph contains only undirected edges, then it has two natural orientations. On the other hand, if $C$ contains some edge which fails to have both S-to-R and R-to-S direction, then $C$ has one natural orientation. Thus there are always either one or two orientations for any cycle. Once an orientation is chosen for a cycle $C$, then each edge (including undirected edges) in $C$ inherits an orientation, which we can call that edge's ``$C$-orientation''. Two cycles $C$ and $D$ are said to have {\bf compatible orientation} if one can choose an orientation for $C$ and an orientation for $D$ such that each edge in their intersection has the same $C$-orientation and $D$-orientation. As shown by example in Figure~\ref{orient}, even two unoriented cycles may have incompatible orientation.

{\bf S-to-R intersection between cycles}. The intersection between two cycles consists of a set of vertex-disjoint components. Two distinct cycles in a DSR graph are said to have S-to-R intersection if they have compatible orientation, and moreover each component of their intersection has odd length. 

{\bf Subgraphs of DSR graphs}. If a DSR graph $G = G_{A, B}$ is associated with a pair $(A, B)$ of $n \times m$ matrices, then given nonempty $\gamma \subset \{1, \ldots, n\}$ and $\delta \subset \{1, \ldots, m\}$, $G(\gamma|\delta)$ will mean $G_{A(\gamma|\delta), B(\gamma|\delta)}$. $G(\gamma|\delta)$ can be regarded as a subgraph of $G$. The definition extends naturally to the case where $A$ and $B$ are replaced with sets of matrices. DSR graphs or subgraphs with an equal number of S- and R-vertices will be referred to as {\bf square}. A square subgraph in which each vertex has exactly one edge incident on it is called a {\bf term subgraph}.

{\bf Associating DSR graphs with functions.} Consider some $f \in \mathcal{D}(X)$, with Jacobian $Df(x)$. Define $\mathcal{G}_{f(x)}$, the set of all DSR graphs associated with $f$ at $x$, as follows:
\[
\mathcal{G}_{f(x)} = \{G_{A, B^T}\,|\,AB = -Df(x)\}\,.
\]
Let the DSR graph $G_{\mathcal{A}, \mathcal{B}}$ be ``associated with $f$'' if $\mathcal{A}, \mathcal{B}$ are sets of matrices of equal dimension such that for each $x \in X$ there exists $A \in \mathcal{A}, B \in \mathcal{B}$ satisfying $-AB^T = Df(x)$. Similarly the DSR graph $G_{\mathcal{A}, \mathcal{B}}$ is associated with $\mathcal{F} \subset \mathcal{D}(X)$ if $\mathcal{A}, \mathcal{B}$ are sets of matrices of equal dimension such that for each $x \in X, f \in \mathcal{F}$ there exists $A \in \mathcal{A}, B \in \mathcal{B}$ satisfying $-AB^T = Df(x)$. Define $\mathcal{G}_{f}$ to be the set of all DSR graphs associated with $f$, and $\mathcal{G}_{\mathcal{F}}$ to be the set of all DSR graphs associated with $\mathcal{F}$.\\

\section{DSR graphs: results}

Define the following conditions on a DSR graph:
\nopagebreak
\begin{quote}
{\bf Condition ($**$)}: It contains no e-cycles. \\
\nopagebreak
{\bf Condition ($*$)}: All e-cycles are s-cycles, and no two e-cycles have S-to-R intersection. 
\end{quote}

Note that Condition ($**$) is more restrictive than Condition ($*$). The key theoretical result underpinning claims in this paper is the following:

\begin{thm}
\label{mainSRgraph}
Given $f \in \mathcal{D}(X)$, assume that there exists $q \in \mathcal{D}_d^{+}(X)$ and $a, b \in X$ ($a \not = b$) such that $f(a)-q(a) = f(b)-q(b)$. Then there exists some $x$ such that each $G\in \mathcal{G}_{f(x)}$ fails Condition~($*$). Thus every $G \in \mathcal{G}_f$ fails Condition~($*$).
\end{thm}
The proof of Theorem~\ref{mainSRgraph} is lengthy and is developed in \cite{banajicraciun2}. It follows from the fact that if some $G \in \mathcal{G}_{f(x)}$ satisfies Condition~($*$), then it can be shown that $-Df(x)$ lies in the closure of the $P$ matrices, and if this is the case at each $x \in X$, then $f-q$ is injective on $X$ for arbitrary $q \in \mathcal{D}_d^{+}(X)$. 
\begin{cor}
\label{mainSRcor}
Given $\mathcal{F} \subset \mathcal{D}(X)$, suppose there exists some $G \in \mathcal{G}_{\mathcal{F}}$ satisfying Condition~($*$). Then $\mathcal{F}^{-}$ is injective. 
\end{cor}
\begin{pf}
This follows from Theorem~\ref{mainSRgraph} because if some $G \in \mathcal{G}_{\mathcal{F}}$ satisfies Condition~($*$), then for each $f \in \mathcal{F}$, there exists some $G \in \mathcal{G}_f$ which satisfies Condition~($*$). \qquad \qed
\end{pf}

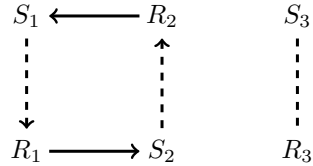
\begin{figure}[h]
\begin{minipage}{\textwidth}
\begin{center}
\begin{tikzpicture}[domain=0:4,scale=0.6]

\node at (1,4) {$S_1$};
\node at (4,4) {$R_2$};
\node at (7,4) {$S_3$};
\node at (1,1) {$R_1$};
\node at (4,1) {$S_2$};
\node at (7,1) {$R_3$};

\draw[<-, line width=0.04cm] (1.5, 4) -- (3.5,4);
\draw[->, line width=0.04cm] (1.5, 1) -- (3.5,1);

\draw[->, dashed, line width=0.04cm] (1, 3.5) -- (1,1.5);
\draw[<-, dashed, line width=0.04cm] (4, 3.5) -- (4,1.5);
\draw[-, dashed, line width=0.04cm] (7, 3.5) -- (7,1.5);

\end{tikzpicture}
\end{center}

\end{minipage}
\caption{\label{notfullynon} The DSR graph shown (edge-labels omitted) is weakly nondegenerate, but not nondegenerate: there is no square subgraph involving S-vertices $S_1$ and $S_3$ which includes both an S-to-R term subgraph and an R-to-S term subgraph.}
\end{figure}

{\bf Extending the results.} It is useful to define the following nondegeneracy conditions on an $n \times m$ DSR graph $G$. $G$ is {\bf weakly nondegenerate} if it contains an $n \times n$ subgraph containing a term subgraph with S-to-R direction and one with R-to-S direction. $G$ is {\bf nondegenerate} if given any subset of the S-vertices, there is a square subgraph including this subset of the S-vertices (and no others) and containing a term subgraph with S-to-R direction and one with R-to-S direction. In other words, given any nonempty $\gamma \subset \{1, \ldots, n\}$, there is a $\delta \subset \{1, \ldots, m\}$ with $|\delta| = |\gamma|$ and such that $G(\gamma|\delta)$ is weakly nondegenerate. Note that in this definition it is important that the {\em same} square subgraph contains a term subgraph with S-to-R direction and one with R-to-S direction (see Figure~\ref{notfullynon}).

\begin{lem}
\label{nondegen}
Consider a pair of $n \times m$ matrices $(A, B)$ and DSR graph $G = G_{A, B}$.\\ 
{\bf 1.} If $AB^T$ is nonsingular then $G$ is weakly nondegenerate.\\
{\bf 2.} Suppose $G$ satisfies Condition~($**$). Then $AB^T$ is nonsingular if and only if $G$ is weakly nondegenerate. 
 
\end{lem}
\begin{pf}
Let $\gamma = \{1, \ldots, n\}$. By the Cauchy-Binet formula, 
\[
\mathrm{det}(AB^T) = \sum_{\substack{\delta \subset \{1, \ldots, m\}\\ |\delta| = n}}A[\gamma|\delta]B[\gamma|\delta]\,.
\]

{\bf 1.} If $AB^T$ is nonsingular, then there exists $\delta \subset \{1, \ldots, m\}$ with $|\delta| = n$ such that $A[\gamma|\delta]B[\gamma|\delta] \not = 0$. So there is at least one nonzero term in $A[\gamma|\delta]$ and similarly in $B[\gamma|\delta]$. But a nonzero term in $A[\gamma|\delta]$ corresponds precisely to an R-to-S term subgraph in $G(\gamma|\delta)$, and similarly a nonzero term in $B[\delta|\gamma]$ corresponds precisely to an S-to-R term subgraph in $G(\gamma|\delta)$. Since $|\gamma| = n$, $G(\gamma|\delta)$ includes all the S-vertices in $G_{A, B}$. 

{\bf 2.} It was proved in part 1 that if $AB^T$ is nonsingular then $G$ is weakly nondegenerate. Since $G$ satisfies Condition~($*$), $A[\gamma|\delta]B[\gamma|\delta] \geq 0$ for each $\gamma \subset \{1, \ldots, n\}$, and $\delta \subset \{1, \ldots, m\}$ with $|\delta| = |\gamma|$. Let $\gamma = \{1, \ldots, n\}$. Since $G$ is weakly nondegenerate, there exists at least one $\delta \subset \{1, \ldots, m\}$ with $|\delta| = n$ such that $G(\gamma|\delta)$ contains a term subgraph with S-to-R direction and one with R-to-S direction. Let these term subgraphs correpond to nonzero terms $T_1$ in the expansion of $A[\gamma|\delta]$ and $T_2$ in the expansion of $B[\gamma|\delta]$. 

Since $G(\gamma|\delta)$ satisfies Condition~($**$), $T_aT_b \geq 0$ for any terms $T_a$ in the expansion of $A[\gamma|\delta]$ and $T_b$ in the expansion of $B[\gamma|\delta]$ (see Lemma 5.1 in \cite{banajicraciun2}), and since $T_1$ and $T_2$ are nonzero, $T_1T_2 > 0$. If $A[\gamma|\delta]=0$, then there must be some term $T_3$ in the expansion of $A[\gamma|\delta]$ such that $T_1T_3 < 0$, contradicting $T_3T_2 \geq 0$. So $A[\gamma|\delta]\not =0$. Similarly $B[\gamma|\delta] \not = 0$. So $A[\gamma|\delta]B[\gamma|\delta] > 0$, i.e. $\mathrm{det}(AB^T) > 0$. \qed

\end{pf}

{\bf Remarks.} As a trivial implication of part 1 of Lemma~\ref{nondegen}, if $G_{A, B}$ contains fewer R-vertices than S-vertices, then $AB^T$ is singular. An application of part 2 of Lemma~\ref{nondegen} is the following: suppose a square matrix $A$ is such that $G_{A, I}$ satisfies Condition ($**$) and is weakly nondegenerate. Then $A$ is nonsingular. Further, for each $B \in \mathcal{Q}(A)$ (i.e. for any $B$ with the same sign pattern as $A$), $G_{B, I}$ is identical, upto edge-labelling, to $G_{A, I}$, and so also satisfies Condition ($**$). In other words, $A$ is sign nonsingular -- that is all matrices with the same sign pattern as $A$ are nonsingular. \\

For a $P$ matrix, every principal minor is nonzero, giving the following result:

\begin{lem}
\label{nondegen1}
Consider a pair of $n \times m$ matrices $(A, B)$ and DSR graph $G = G_{A, B}$.\\ 
{\bf 1.} If $AB^T$ is a $P$ matrix, then $G$ is nondegenerate. \\
{\bf 2.} Suppose $G$ satisfies Condition~($**$). Then $AB^T$ is a $P$ matrix if and only if $G$ is nondegenerate. 
\end{lem}
\begin{pf}
Let $\gamma \subset \{1, \ldots, n\}$ be nonempty. By the Cauchy-Binet formula, 
\[
(AB^T)[\gamma] = \sum_{\substack{\delta \subset \{1, \ldots, m\}\\ |\delta| = |\gamma|}}A[\gamma|\delta]B[\gamma|\delta]\,,
\]

{\bf 1.} If $AB^T$ is a $P$ matrix, then for each such $\gamma$, there exists $\delta \subset \{1, \ldots, m\}$ with $|\delta| = |\gamma|$ such that $A[\gamma|\delta]B[\gamma|\delta] \not = 0$. So there is at least one nonzero term in $A[\gamma|\delta]$ and similarly in $B[\gamma|\delta]$. But a nonzero term in $A[\gamma|\delta]$ corresponds precisely to an R-to-S term subgraph in $G(\gamma|\delta)$, and similarly a nonzero term in $B[\delta|\gamma]$ corresponds precisely to an S-to-R term subgraph in $G(\gamma|\delta)$. By definition $G(\gamma|\delta)$ includes all the S-vertices indexed by $\gamma$ and no others, so $G$ is nondegenerate.

{\bf 2.} It was proved in part 1 that if $AB^T$ is a $P$ matrix, then $G$ is nondegenerate. Since $G$ satisfies Condition~($*$), $A[\gamma|\delta]B[\gamma|\delta] \geq 0$ for each nonempty $\gamma \subset \{1, \ldots, n\}$, and $\delta \subset \{1, \ldots, m\}$ with $|\delta| = |\gamma|$. Let $\gamma \subset \{1, \ldots, n\}$ be nonempty. Since $G$ is nondegenerate, there exists at least one $\delta \subset \{1, \ldots, m\}$ with $|\delta| = |\gamma|$ such that $G(\gamma|\delta)$ contains a term subgraph with S-to-R direction and one with R-to-S direction. Applying the arguments in part 2 of Lemma~\ref{nondegen} to $A[\gamma|\delta]$ and $B[\gamma|\delta]$, gives $A[\gamma|\delta]B[\gamma|\delta] > 0$, i.e. $(AB^T)[\gamma] > 0$. Since $\gamma$ was arbitrary, $AB^T$ is a $P$ matrix. \qed
\end{pf}

{\bf Remark.} An application is the following. Suppose a square matrix $A$ is such that $G_{A, I}$ satisfies Condition ($**$) and is nondegenerate. Then $A$ is a $P$ matrix. Further, for each $B \in \mathcal{Q}(A)$, $G_{B, I}$ is identical, upto edge-labelling, to $G_{A, I}$, and so also satisfies Condition ($**$) and is nondegenerate. In other words, all matrices in $\mathcal{Q}(A)$ are $P$ matrices.

\begin{thm}
\label{nond}
Consider some $f \in \mathcal{D}(X)$.\\
{\bf 1.} If, at each $x \in X$, there exists a nondegenerate DSR graph $G \in \mathcal{G}_{f(x)}$ satisfying Condition~($**$), then $f^{-}$ and $f$ are injective.\\
{\bf 2.} If $X$ is open and, at each $x \in X$, there exists a weakly nondegenerate DSR graph $G \in \mathcal{G}_{f(x)}$ satisfying Condition~($**$), then $f^{-}$ and $f$ are injective.

\end{thm}
\begin{pf}
In both cases, since Condition~($**$) implies Condition~($*$), by Theorem~\ref{mainSRgraph}, $f^{-}$ is injective.\\ 
{\bf 1.} By Lemma~\ref{nondegen1}, $-Df(x)$ is a $P$ matrix at each $x \in X$. Thus $-f$, and hence $f$, is injective on $X$ \cite{gale}. \\
{\bf 2.} By Lemma~\ref{nondegen}, $-Df(x)$, and hence $Df(x)$, is nonsingular at each $x \in X$. Further, as $G$ certainly satisfies Condition~($*$), $-Df(x)$ is in the closure of the $P$ matrices at each $x \in X$. Thus, since $X$ is open, $f$ is injective on $X$ (see Theorem~4w in \cite{gale} and Appendix B in \cite{banajicraciun2}). \qed
\end{pf}

\section{The Jacobian DSR graph: relationship between I-graph and DSR graph results}

Given a function $f \in \mathcal{D}(X)$, writing $f = f \circ \mathrm{id}$ (where $\mathrm{id}$ is the identity on $X$), gives a natural factorisation of the Jacobian at each point $Df(x) = Df(x)\,I$, leading to DSR graphs $G_{Df(x), -I}$. Given any square matrix $M$, the particular DSR graph $G_{M, -I}$ will be termed the {\bf Jacobian DSR graph} corresponding to $M$, or JDSR graph for short. Note that JDSR graphs are always square.

\begin{thm}
\label{mainimp}
Consider a square matrix $M$ with corresponding I-graph $H = H_M$ and JDSR graph $G = G_{M, -I}$. The following statements are equivalent:
\begin{enumerate}
\item $H$ contains a positive (resp. negative) cycle.
\item $G$ contains an e-cycle (resp. o-cycle). 
\end{enumerate}
\end{thm}

\begin{cor}
\label{mainimplication}
Consider some $f \in \mathcal{D}(X)$. At each point $x \in X$, associate with $f$ the I-graph $H_{Df(x)}$ and the JDSR graph $G_{Df(x), -I}$. Then $H_{Df(x)}$ contains no positive cycles if and only if $G_{Df(x), -I}$ satisfies Condition~($**$).
\end{cor}
\begin{pf}
This follows trivially from Theorem~\ref{mainimp}. \qed
\end{pf}

\begin{pot1}
The equivalence between existence of a positive cycle in $H$ and an e-cycle in $G$ will be proved. The equivalence between existence of a negative cycle in $H$ and an o-cycle in $G$ follows analogously. A directed edge from vertex $j$ to vertex $i$ in $H$ will be termed $h_{ij}$. Similarly $g_{ij}$ will refer to an edge in $G$ between S-vertex $i$ and R-vertex $j$ with arrows/lines above indicating direction. 

{\bf Statement 1 implies statement 2.} Assume the existence of a positive $n$-cycle ($n \geq 2$) in $H$:
\[
C_H = \{h_{i_1i_2}, h_{i_2i_3},\ldots,h_{i_ni_1}\}\,,
\]
where $i_j \not = i_k$ for $j \not = k$. Let $i_{k+1}$ mean $i_{(k \mod n) + 1}$. Since $C_H$ is positive,
\[
\mathrm{sign}(C_H) = \prod_{k=1}^n\mathrm{sign}(h_{i_ki_{k+1}}) = 1.
\]

An edge $h_{i_ki_{k+1}}$ corresponds to an entry $M_{i_ki_{k+1}}$ in $M$, and hence to an edge $\overleftarrow{g}_{i_ki_{k+1}}$ in $G$. Since $i_j \not = i_k$ for $j \not = k$, no two of these edges share a vertex. Moreover $\mathrm{sign}(h_{i_ki_{k+1}}) = \mathrm{sign}(\overleftarrow{g}_{i_ki_{k+1}})$, so
\[
\prod_{k=1}^n\mathrm{sign}(\overleftarrow{g}_{i_ki_{k+1}}) = 1.
\]
Now the JDSR graph, by definition contains negative edges $\overrightarrow{g}_{i_k, i_k}$. Thus there is the following cycle of length $2n$ in $G$:
\[
C_G = \{\overrightarrow{g}_{i_1i_1}, \overleftarrow{g}_{i_1i_2}, \overrightarrow{g}_{i_2i_2}, \overleftarrow{g}_{i_2i_3},\ldots,\overrightarrow{g}_{i_ni_n}, \overleftarrow{g}_{i_ni_1}\}\,.
\]
So 
\[
\mathrm{sign}(C_G) = \left(\prod_{k=1}^n\mathrm{sign}(\overleftarrow{g}_{i_ki_{k+1}})\right)\left(\prod_{k =1}^n\mathrm{sign}(\overrightarrow{g}_{i_ki_k})\right) = (-1)^{n}.
\]
Since $|C_G|/2 = n$, the parity of $C_G$ is 
\[
P(C_G) = (-1)^{|C_{G}|/2}\mathrm{sign}(C_G) = (-1)^{n}(-1)^{n} = 1
\]
and thus $C_G$ is an e-cycle.

{\bf Statement 2 implies statement 1.} Assume the existence of an e-cycle $C_G$ of length $2n$ in $G$. Since the only edges in $G$ with S-to-R direction are edges of the form $\overrightarrow{g}_{kk}$, such an e-cycle must take the form
\[
C_G = \{\overrightarrow{g}_{i_1i_1}, \overleftarrow{g}_{i_1i_2}, \overrightarrow{g}_{i_2i_2}, \overleftarrow{g}_{i_2i_3},\ldots,\overrightarrow{g}_{i_ni_n}, \overleftarrow{g}_{i_ni_1}\}\,,
\]
for some set of indices $K = \{i_1, \ldots, i_n\}$. As before, by the definition of a cycle, $i_j \not = i_k$ for $j \not = k$. As before, $P(C_G) = 1$ implies $\mathrm{sign}(C_G) = (-1)^{n}$. But the fact that edges $\overrightarrow{g}_{i_ki_k}$ are negative means that $\prod_{k =1}^n\mathrm{sign}(\overrightarrow{g}_{i_ki_k}) = (-1)^{n}$. So 
\[
\prod_{k =1}^n\mathrm{sign}(\overleftarrow{g}_{i_ki_{k+1}}) = 1.
\]
The existence of edges $\overleftarrow{g}_{i_ki_{k+1}}$ in $C_G$ implies the existence of the $n$-cycle in $H$:
\[
C_H = \{h_{i_1i_2}, h_{i_2i_3},\ldots,h_{i_ni_1}\}\,.
\]
Since $\prod_{k =1}^n\mathrm{sign}(\overleftarrow{g}_{i_ki_{k+1}}) = 1$, this implies that $\prod_{k =1}^n\mathrm{sign}(h_{i_ki_{k+1}}) = 1$. Thus $C_H$ is a positive cycle. \qquad \qed
\end{pot1}

The theorem tells us that any conclusions that can be drawn from the absence or presence of positive (resp. negative) cycles in $H_{\mathcal{F}}$, can also be derived from the JDSR graph. It will be shown by example that the converse is not true: for example, there are systems with JDSR graph which satisfy Condition~($*$), but which have positive cycles in $H_{\mathcal{F}}$. In fact, defining:
\begin{enumerate}
\item[C1.] Functions whose I-graphs have no positive cycles;
\item[C2.] Functions whose JDSR graphs satisfy Condition~($*$);
\item[C3.] Functions for which there exists a DSR graph which satisfies Condition~($*$),
\end{enumerate}
then C1 is a proper subset of C2, and C2 is a proper subset of C3.\\ 

{\bf Proofs of I-graph results.} Theorems~\ref{mainIgraph}~and~\ref{mainIgrapha} become immediate corollaries of Theorem~\ref{mainSRgraph} and Corollary~\ref{mainimplication}:


\begin{pot2}
Since $f - q$ is noninjective for some $q \in \mathcal{D}_d^{+}(X)$, by Theorem~\ref{mainSRgraph}, there exists $c \in X$ such that $G_{Df(c), -I}$ (and indeed any other DSR graph $G \in \mathcal{G}_{f(c)}$) fails Condition~($*$). By Corollary~\ref{mainimplication}, the I-graph $H_{Df(c)}$ (and hence $H_f$) contains a positive cycle. \qed
\end{pot2}

\begin{pot3}
If the Jacobian $Df(x)$ has negative diagonal elements, then the JDSR graph $G = G_{Df(x), -I}$ contains S-to-R and R-to-S term subgraphs involving precisely edges of the form $S_i \!-\! R_i$. Given any nonempty $\gamma \subset \{1, \ldots, n\}$, $G(\gamma|\gamma)$ is thus weakly nondegenerate, and so $G$ is nondegenerate. From Corollary~\ref{mainimplication} if $H_{Df(x)}$ has no positive cycles for any $x \in X$, then $G$ satisfies Condition~($**$). Thus, by Theorem~\ref{nond}, $f$ is injective on $X$. \qed
\end{pot3}

\section{Examples}

{\bf Example 1. Choosing a factorisation.} Define $x = [x_1, x_2]^T \in \mathbb{R}^2$ and let $X$ be any rectangular subset of $\mathbb{R}^2$. Let $f_1(x_1)$ and $f_2(x_2)$ be real functions such that $f_1^{'}(x_1) = \frac{\mathrm{d}f_1}{\mathrm{d}x_1} > 0$, $f_2^{'}(x_2) = \frac{\mathrm{d}f_2}{\mathrm{d}x_2} > 0$ for all $x \in X$. Consider the function $f:X \to \mathbb{R}^2$ defined by 
\begin{equation}
\label{ex1}
f(x) = \left[\begin{array}{c}-f_1(x_1) - f_2(x_2)/2\\-f_1(x_1) - f_2(x_2)
\end{array}\right],
\end{equation}
with Jacobian
\[
Df(x) = \left[\begin{array}{rr}-f_1^{'}(x_1) & -f_2^{'}(x_2)/2\\-f_1^{'}(x_1) & -f_2^{'}(x_2)\end{array}\right].
\]

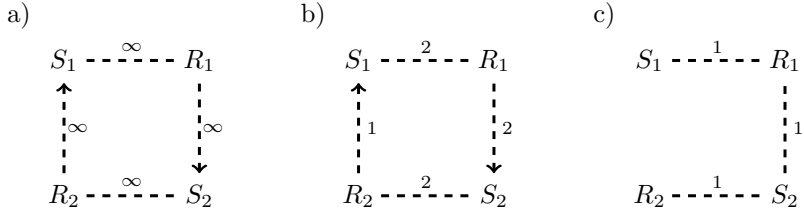
\begin{figure}[h]
\begin{center}
\begin{minipage}{0.31\textwidth}
\begin{center}
\begin{tikzpicture}[domain=0:4,scale=0.6]

\node at (0,5) {$\mathrm{a)}$};

\node at (1,4) {$S_1$};
\node at (4,4) {$R_1$};
\node at (4,1) {$S_2$};
\node at (1,1) {$R_2$};

\draw[-, dashed, line width=0.04cm] (1.5,4.0) -- (3.5, 4.0);
\draw[-, dashed, line width=0.04cm] (1.5,1.0) -- (3.5, 1.0);

\draw[->, dashed, line width=0.04cm] (1.0, 1.5) -- (1.0, 3.5);
\draw[<-, dashed, line width=0.04cm] (4.0, 1.5) -- (4.0, 3.5);

\node at (2.5, 4.3) {$\scriptstyle{\infty}$};
\node at (2.5, 1.3) {$\scriptstyle{\infty}$};

\node at (1.3, 2.5) {$\scriptstyle{\infty}$};
\node at (4.3, 2.5) {$\scriptstyle{\infty}$};

\end{tikzpicture}
\end{center}
\end{minipage}
\begin{minipage}{0.31\textwidth}
\begin{center}
\begin{tikzpicture}[domain=0:4,scale=0.6]

\node at (0,5) {$\mathrm{b)}$};

\node at (1,4) {$S_1$};
\node at (4,4) {$R_1$};
\node at (4,1) {$S_2$};
\node at (1,1) {$R_2$};

\draw[-, dashed, line width=0.04cm] (1.5,4.0) -- (3.5, 4.0);
\draw[-, dashed, line width=0.04cm] (1.5,1.0) -- (3.5, 1.0);

\draw[->, dashed, line width=0.04cm] (1.0, 1.5) -- (1.0, 3.5);
\draw[<-, dashed, line width=0.04cm] (4.0, 1.5) -- (4.0, 3.5);

\node at (2.5, 4.3) {$\scriptstyle{2}$};
\node at (2.5, 1.3) {$\scriptstyle{2}$};

\node at (1.3, 2.5) {$\scriptstyle{1}$};
\node at (4.3, 2.5) {$\scriptstyle{2}$};

\end{tikzpicture}
\end{center}
\end{minipage}
\begin{minipage}{0.31\textwidth}
\begin{center}
\begin{tikzpicture}[domain=0:4,scale=0.6]

\node at (0,5) {$\mathrm{c)}$};

\node at (1,4) {$S_1$};
\node at (4,4) {$R_1$};
\node at (4,1) {$S_2$};
\node at (1,1) {$R_2$};

\draw[-, dashed, line width=0.04cm] (1.5,4.0) -- (3.5, 4.0);
\draw[-, dashed, line width=0.04cm] (1.5,1.0) -- (3.5, 1.0);

\draw[-, dashed, line width=0.04cm] (4.0, 1.5) -- (4.0, 3.5);

\node at (2.5, 4.3) {$\scriptstyle{1}$};
\node at (2.5, 1.3) {$\scriptstyle{1}$};

\node at (4.3, 2.5) {$\scriptstyle{1}$};

\end{tikzpicture}
\end{center}
\end{minipage}
\end{center}
\caption{\label{decomp} The DSR graphs corresponding to three factorisations (Eqs.~\ref{decomp1},~\ref{decomp2}~and~\ref{decomp3}) of the Jacobian of the function $f$ in Eq.~\ref{ex1}. DSR graphs a) and b) fail Condition~($*$). DSR graph c) is in fact a tree and satisfies Condition~($**$). It is also nondegenerate.}
\end{figure}

Consider three factorisations of $Df(x)$: 
\begin{eqnarray}
\label{decomp1}
Df(x) & =& \left[\begin{array}{rr}-f_1^{'}(x_1) & -f_2^{'}(x_2)/2\\-f_1^{'}(x_1) & -f_2^{'}(x_2)\end{array}\right]\left[\begin{array}{rr}1 & 0 \\0 & 1\end{array}\right],\\
\label{decomp2}
Df(x) & = & \left[\begin{array}{rr}-2 & -1\\-2 & -2\end{array}\right]\left[\begin{array}{cc}f_1^{'}(x_1)/2 & 0 \\0& f_2^{'}(x_2)/2\end{array}\right],\\
\label{decomp3}
Df(x) & = & \left[\begin{array}{rr}-1 & 0\\-1 &-1\end{array}\right]\left[\begin{array}{cc}f_1^{'}(x_1) & f_2^{'}(x_2)/2\\0 & f_2^{'}(x_2)/2\end{array}\right].
\end{eqnarray}
These three factorisations give the three DSR graphs shown in Figure~\ref{decomp}. The first two fail Condition~($*$), while the third satisfies Condition~($**$) (indeed it is a tree). Moreover, by inspection it is nondegenerate. Thus, by Theorem~\ref{nond}, $f^{-}$ and $f$ are injective on $X$. It is not obvious {\em a priori}, that the third factorisation is likely to be the most useful.

\begin{figure}[h]
\begin{center}
\begin{tikzpicture}[domain=0:4,scale=0.65]

\node at (0,3.8) {$x_3$};
\node at (3,3.8) {$x_2$};
\node at (1.5,1.2) {$x_1$};

\draw[->, line width=0.04cm] (2.5,4.0) .. controls (1.9,4.3) and (1.1,4.3) .. (0.5,4.0);
\draw[<-, line width=0.04cm] (2.5,3.7) .. controls (1.9,3.4) and (1.1,3.4) .. (0.5,3.7);

\draw[->, dashed, line width=0.04cm] (2.8,3.5) -- (1.7, 1.5);
\draw[->, line width=0.04cm] (1.3, 1.5) -- (0.2, 3.5);

\node at (7,4) {$R_1$};
\node at (10,4) {$x_3$};
\node at (13,4) {$R_3$};
\node at (7,1) {$x_1$};
\node at (10,1) {$R_2$};
\node at (13,1) {$x_2$};

\draw[->, line width=0.04cm] (7.5,4.0) -- (9.5, 4.0);
\draw[-, dashed, line width=0.04cm] (7.0,1.5) -- (7.0,3.5);
\draw[<-, dashed, line width=0.04cm] (7.5,1.0) -- (9.5, 1.0);

\draw[-, dashed, line width=0.04cm] (10.5,4.0) -- (12.5, 4.0);

\draw[-, dashed, line width=0.04cm] (10.5,1.0) -- (12.5, 1.0);

\draw[->, line width=0.04cm] (10.0,1.5) -- (10.0,3.5);
\draw[<-, line width=0.04cm] (13.0,1.5) -- (13.0,3.5);

\node at (11.5, 4.3) {$\scriptstyle{2}$};
\node at (11.5, 0.7) {$\scriptstyle{1}$};
\node at (9.7, 2.5) {$\scriptstyle{1}$};
\node at (13.3, 2.5) {$\scriptstyle{2}$};

\node at (8.5, 4.3) {$\scriptstyle{\infty}$};
\node at (8.5, 0.7) {$\scriptstyle{\infty}$};
\node at (6.7, 2.5) {$\scriptstyle{\infty}$};

\end{tikzpicture}
\end{center}
\caption{\label{pic2}{\em Left.} The I-graph $H_{Df(x)}$ associated with the function $f$ in Example 2 at any $x$. Negative self-edges have been omitted. $H_{Df(x)}$ contains a positive cycle, and so Theorem~\ref{mainIgraph} cannot be used to draw any conclsions. {\em Right.} The JDSR graph $G_{Df(x), -I}$ at each $x$ satisfies Condition~($*$) and so $f^{-}$ is injective.}
\end{figure}
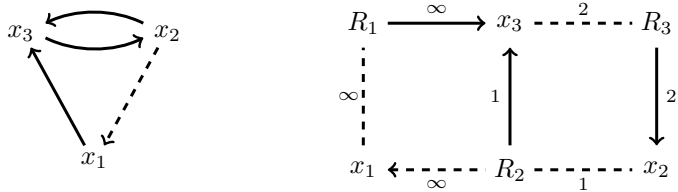

{\bf Example 2. Functions with some linear terms.} Define $x = [x_1, x_2, x_3]^T$ and $X$ to be some rectangular subset of ${R}^3$. Consider the function $f:X \to \mathbb{R}^3$ defined by $F = [f_1(x_1, x_2), x_3-x_2, f_2(x_1)+2(x_2-x_3)]^T$ where $\frac{\partial f_1}{\partial x_1} < 0$, $\frac{\partial f_1}{\partial x_2} < 0$ and $\frac{\partial f_2}{\partial x_1} > 0$. The system has Jacobian with structure
\[
Df(x) = \left[ \begin{array}{rrr}-a &  -b &  0\\0 & -1 &  1\\c &  2 &  -2\end{array}\right]
\]
where $a,b,c> 0$. The system has I-graph and JDSR graph shown in Figure~\ref{pic2}. The JDSR graph satisfies Condition~($*$), and so $f^{-}$ is injective on $X$. On the other hand the I-graph contains a positive cycle, and cannot be directly used to draw this conclusion. This example illustrates that even only using the JDSR graph can give stronger results than using the I-graph alone, as edge-labels in the JDSR graph provide information not in the I-graph.

\section{Conclusions}

A number of graph-theoretic techniques for deciding on injectivity of a set of functions have been described and applied to examples. It has been shown that examining a particular DSR graph, termed the JDSR graph, allows stronger conclusions about injectivity than are possible from the I-graph alone. Note that Theorem~\ref{mainimp} implies that {\em any} conclusions which follow from the presence or absence of cycles in an I-graph, can equally be drawn from DSR graphs. Applications of DSR graph techniques to questions going beyond injectivity will also be explored in future work.

A theoretical difficulty is that there is no unique way of associating DSR graphs with functions. Thus an important challenge is to find systematic ways -- either analytical or algorithmic -- of choosing factorisations of Jacobians, and hence DSR graphs, to allow the strongest conclusions about injectivity. In many contexts, natural structures can be exploited in associating DSR graphs with dynamical systems. This was illustrated via a number of nontrivial examples drawn from the applied literature in \cite{banajicraciun2}. Further real examples will be presented in forthcoming work.

\section*{References}

\bibliographystyle{unsrt}

\appendix

\section{Dual results}
\label{dualapp}

The following results are collected for completeness. Their proofs follow closely the corresponding results in brackets and are omitted.

\begin{thm}[corresponding to Theorem~\ref{mainIgraph}]
\label{mainIgraphdual}
Given $f \in \mathcal{D}(X)$, suppose there exists some $q \in \mathcal{D}_d^{+}(X)$, and $a, b \in X$  ($a \not = b$) such that $f(a)+q(a) = f(b)+q(b)$. Then there exists $c \in X$ such that $H_{-Df(c)}$ contains a positive cycle, and thus $H_{-f}$ contains a positive cycle.
\end{thm}

\begin{thm}[corresponding to Theorem~\ref{mainIgrapha}]
\label{mainIgraphadual}
Given $f \in \mathcal{D}(X)$ such that the Jacobian $Df$ has positive diagonal elements (i.e. $\frac{\partial f_i}{\partial x_i} > 0$  at each $x \in X$), suppose there exist $a, b \in X$ ($a \not = b$) such that $f(a)= f(b)$.  Then there exists $c \in X$ such that $H_{-Df(c)}$ contains a positive cycle, and thus $H_{-f}$ contains a positive cycle.
\end{thm}

\begin{cor}[corresponding to Corollary~\ref{corIgraph}]
\label{corIgraphdual}
For some $\mathcal{F} \subset \mathcal{D}(X)$, assume that $H_\mathcal{-F}$ contains no positive cycles. \\
{\bf 1.} Then $\mathcal{F}^{+}$ is injective. \\
{\bf 2.} Assume in addition that $\frac{\partial f_i}{\partial x_i} > 0$ for each $f \in \mathcal{F}, x \in X$. Then $\mathcal{F}^{+} \cup \mathcal{F}$ is injective.
\end{cor}

\begin{thm}[corresponding to Theorem~\ref{mainSRgraph}]
\label{mainSRgraphdual}
Given $f \in \mathcal{D}(X)$, assume that there exists $q \in \mathcal{D}_d^{+}(X)$ and $a, b \in X$ ($a \not = b$) such that $f(a)+q(a) = f(b)+q(b)$. Then there exists some $x$ such that each $G\in \mathcal{G}_{-f(x)}$ fails Condition~($*$). Thus every $G \in \mathcal{G}_{-f}$ fails Condition~($*$).
\end{thm}

\begin{cor}[corresponding to Corollary~\ref{mainSRcor}]
\label{mainSRcordual}
Given $\mathcal{F} \subset \mathcal{D}(X)$, suppose there exists some $G \in \mathcal{G}_{\mathcal{-F}}$ satisfying Condition~($*$). Then $\mathcal{F}^{+}$ is injective. 
\end{cor}

\begin{thm}[corresponding to Theorem~\ref{nond}]
\label{nonddual}
Consider some $f \in \mathcal{D}(X)$.\\
{\bf 1.} If, at each $x \in X$, there exists a nondegenerate DSR graph $G \in \mathcal{G}_{-f(x)}$ satisfying Condition~($**$), then $f^{+}$ and $f$ are injective.\\
{\bf 2.} If $X$ is open and, at each $x \in X$, there exists a weakly nondegenerate DSR graph $G \in \mathcal{G}_{-f(x)}$ satisfying Condition~($**$), then $f^{+}$ and $f$ are injective.
\end{thm}

\begin{cor}[corresponding to Corollary~\ref{mainimplication}]
\label{mainimplicationdual}
Consider some $f \in \mathcal{D}(X)$. At each point $x \in X$, consider the I-graph $H_{-Df(x)}$ and the JDSR graph $G_{Df(x), I}$. Then $H_{-Df(x)}$ contains no positive cycles if and only if $G_{Df(x), I}$ satisfies Condition~($**$).
\end{cor}

\end{document}